\documentclass{elsart}

\usepackage{graphicx}
\usepackage{amssymb,amsmath,amsfonts,latexsym,euscript}

\usepackage{amssymb}

\def\d{{\rm d}}
\def\E{{\rm e}}
\def\i{{\rm i}}

\newtheorem{theorem}{Theorem}

\newtheorem{problem}{Problem}
\newtheorem{corollary}{Corollary}

\begin{document}

\begin{frontmatter}



\title{On an optimal quadrature formula for approximation of Fourier integrals in the space $L_2^{(1)}$}


\author[label1,label2]{Abdullo R. Hayotov},
\ead{hayotov@mail.ru}
\author[label1]{Soomin Jeon},
\ead{soominjeon@kaist.ac.kr}
\author[label1]{Chang-Ock Lee}
\ead{colee@kaist.edu}

\address[label1]{Department of Mathematical Sciences, KAIST, 291 Daehak-ro, Yuseong-gu, Daejeon 34141, Republic of Korea}
\address[label2]{V.I.Romanovskiy Institute of Mathematics, Uzbekistan Academy of Sciences,
81, M.Ulugbek str., Tashkent 100170, Uzbekistan}

\begin{abstract}
This paper deals with the construction of an optimal quadrature formula for the approximation of Fourier integrals
in the Sobolev space $L_2^{(1)}[a,b]$ of non-periodic, complex valued functions which are square integrable with first order derivative.
Here the quadrature sum consists of linear combination of the given function values in a uniform grid. The difference between the integral and the quadrature sum is estimated by the norm of the error functional. The optimal quadrature formula is obtained by minimizing the norm of the error functional
with respect to coefficients. Analytic formulas for optimal coefficients can also be obtained using discrete analogue of the differential operator $\d^2/\d x^2$. In addition, the convergence order of the optimal quadrature formula is studied.
It is proved that the obtained formula is exact for all linear polynomials. Thus,
it is shown that the convergence order of the optimal quadrature formula for functions of the space $C^2[a,b]$ is $O(h^2)$.
Moreover, several numerical results are presented and the obtained optimal quadrature formula
is applied to reconstruct the X-ray Computed Tomography image by approximating Fourier transforms.

\end{abstract}

\end{frontmatter}

\section{Introduction}

In practice, since we have discrete values of an integrand, the Fourier transforms are reduced
to an approximation of the integral of type
\begin{equation}\label{FCoef}
I(\varphi)=\int\limits_0^1\E^{2\pi\i\omega x}\varphi(x)\d x
\end{equation}
with $\omega\in \mathbb{R}$. For example, the problem of X-ray Computed Tomography (CT) is to reconstruct
the function from its Radon transform. One of the widely used analytic methods in CT image reconstruction is the filtered back-projection method
in which the Fourier transforms are used (see \cite[Chapter 3]{KakSlaney88} or formulas (\ref{eq4.8_1})-(\ref{eq4.8_3}) of section 4.2).

It should be recalled that integrals of type (\ref{FCoef}) with strongly oscillating integrands
are used in applications of mathematics and other sciences. They are mainly calculated using special effective methods
of numerical integration (for review see, for example, \cite{AvdMal89,BabVitPrag69,BakhVas68,Filon28,IserNor05,Mil98,MilStan14,NovUllWoz15,Olver08,Shad99,XuMilXiang},
and references therein).

Based on Sobolev's method, the problem of the construction of optimal quadrature formulas
for numerical calculation of Fourier coefficients (\ref{FCoef})
with $\omega\in \mathbb{Z}$ in Hilbert spaces $L_2^{(m)}$ and $W_2^{(m,m-1)}$ was studied in \cite{BolHayShad17} and \cite{BolHayMilShad17}, respectively. In these works, explicit formulas of optimal coefficients were obtained for $m\geq 1$.
In particular, for $m=1$, the convergence order of optimal quadrature formulas was studied.

Recently, in \cite{ZhangNovak19} the optimal quadrature formulas were studied for integrals with arbitrary weights in Sobolev space $H^1([0,1])$.
General formulas were obtained for the worst-case error depending on nodes. Especially, when calculating Fourier coefficients of the form
(\ref{FCoef}) with real $\omega$, it was proved that
equidistant nodes are optimal if $n\geq 2.7|\omega|+1$, where $n$ is the number of nodes in the quadrature formula.

It should be noted that for numerical calculation of the integral (\ref{FCoef}) with real $\omega$, a quadrature formula with explicit coefficients
is needed. Therefore, in this paper, we study the construction of optimal quadrature formulas in the sense of Sard
for the approximation of Fourier integrals of the form (\ref{FCoef}) with $\omega\in \mathbb{R}$ in the Sobolev space
of non-periodic square integrable functions with the first order derivative.
We obtain explicit formulas for optimal coefficients and calculate the norm of the error functional of the optimal quadrature formula.
We note that the obtained optimal quadrature formula can be used to approximate Fourier integrals
and reconstruct a function from its discrete Radon transform.

The rest of the paper is organized as follows. In Section 2, an optimal quadrature formula in the sense of Sard is constructed to approximate Fourier integrals in the space $L_2^{(1)}[0,1]$. In Section 3, the results of Section 2 are extended to the case of arbitrary interval $[a,b]$ by linear transformation.
That is, an optimal quadrature formula is obtained for approximate Fourier integrals in the space $L_2^{(1)}[a,b]$.
Finally, in Section 4 the obtained quadrature formula is applied to the approximation of Fourier transforms of a function using the given values of
the function and to the reconstruction of the X-ray CT image.

\section{Construction of optimal quadrature formula for the interval $[0,1]$}
\setcounter{equation}{0}

Consider the quadrature formula
\begin{equation}\label{eq1}
\int\limits_0^1\E^{2\pi\i\omega x}\varphi(x)\d x\cong \sum\limits_{\beta=0}^NC_{\beta}\varphi(h\beta)
\end{equation}
with the error
\begin{equation}
\label{eq2}
(\ell,\varphi)=\int\limits_0^1\E^{2\pi \omega \i x}\varphi(x)\d x -\sum\limits_{\beta  = 0}^N {C_\beta  \varphi (h\beta)},
\end{equation}
where
$$
(\ell,\varphi)=\int\limits_{ - \infty }^\infty  {\ell (x)\varphi(x)\d x},
$$
and the corresponding error functional
\begin{equation}\label{eq3}
\ell (x) =\E^{2\pi \i\omega x}\varepsilon _{[0,1]}(x)-\sum\limits_{\beta  = 0}^N
{C_\beta} \delta (x - h\beta).
\end{equation}
Here, $C_{\beta}$ are coefficients of the formula (\ref{eq1}), $h=1/N$, $N\in \mathbb{N}$, $\i^2=-1$,
$\omega\in \mathbb{R}$ with $\omega\neq 0$, $\varepsilon_{[0,1]}(x)$ is the characteristic function of the interval $[0,1]$, and $\delta$
is  the Dirac's delta-function. The function $\varphi$ belongs to the Sobolev space $L_2^{(1)}[a,b]$ of complex valued functions which are
defined in the interval $[a,b]$ and square integrable with the first order derivative. In this space, the inner product is defined as
\begin{equation}\label{eq4}
\langle\varphi,\psi\rangle =\int\limits_a^b\varphi'(x)\bar\psi'(x)\d x,
\end{equation}
where $\bar\psi$ is the complex conjugate function for the function $\psi$ and
the norm of the function $\varphi$ is denoted by
$$
\|\varphi\|_{L_2^{(1)}[a,b]}=\langle\varphi,\varphi\rangle^{1/2}.
$$

We note that the coefficients $C_{\beta}$ in the formula (\ref{eq1}) vary by $\omega$ and $h$, that is $C_{\beta}=C_{\beta}(\omega,h)$.

The error (\ref{eq2}) in the quadrature formula (\ref{eq1}) is a linear functional
in $L_2^{(1)*}[0,1]$, where $L_2^{(1)*}[0,1]$ is the
conjugate space for the space $L_2^{(1)}[0,1]$.

The absolute value of the error (\ref{eq2}) is estimated by Cauchy-Schwarz inequality  as
$$
|(\ell,\varphi)|\leq \|\varphi\|_{L_2^{(1)}[0,1]}\cdot
\|\ell\|_{L_2^{(1)*}[0,1]},
$$
where
\begin{equation}\label{eq5}
\left\| \ell \right\|_{L_2^{(1)*}[0,1]} = \mathop {\sup
}\limits_{\left\| {\varphi } \right\|_{L_2^{(1)}[0,1]} = 1}
\left| {\left( {\ell,\varphi } \right)} \right|
\end{equation}
is the norm of the error functional (\ref{eq3}).

In the sense of Sard \cite{Sard}, the problem of construction of the optimal quadrature formula (\ref{eq1})
is to find the minimum of the norm (\ref{eq5}) of the error
functional $\ell$ by coefficients $C_{\beta}$ when nodes
are fixed. Here, we note that distances between adjacent nodes in the formula (\ref{eq1}) are the same.
For the quadrature formulas of the form (\ref{eq1}) with $\omega=0$, this problem was first studied by Sard
in $L_2^{(m)}$ space for some $m$, where $L_2^{(m)}$ is the space of real-valued functions which are square integrable with $m$th generalized derivative.
Also this problem for the case $\omega=0$ has been investigated by many authors using splines, $\phi-$function and Sobolev methods. For example, see
\cite{IBab,CatCom,GhOs,FLan,ShadHay11,Sobolev06,Sobolev74,SobVas} and references therein.

Therefore, in order to construct optimal quadrature formulas of the form (\ref{eq1}) in the sense of Sard
in the space $L_2^{(1)}[0,1]$, the following problem needs to be solved.

\begin{problem}\label{Prob1}
Find the coefficients $\mathring{C}_\beta $ that
satisfy the equality
\begin{equation}\label{eq6}
\left\| {\mathring{\ell}}\right\|_{L_2^{(1)*}[0,1]} = \mathop {\inf }\limits_{C_{_\beta  } } \left\| \ell \right\|_{L_2^{(1)*}[0,1]}.
\end{equation}
\end{problem}

In this section we solve Problem \ref{Prob1} for the case $\omega\in \mathbb{R}$ with  $\omega\neq 0$ by finding
the norm (\ref{eq5}) and minimizing it by coefficients $C_{\beta}$.

\subsection{The norm of the error functional (\ref{eq3})}

To find the norm (\ref{eq5}), we use \emph{the extremal function} $\psi_{\ell}$ for the error functional $\ell$ (see \cite{Sobolev74,SobVas})
that satisfies the following equality:
\begin{equation}\label{eq7}
 \left({\ell,\psi _\ell} \right) = \left\| \ell \right\|_{L_2^{(1)*}[0,1]} \cdot
 \left\|\psi_\ell \right\|_{L_2^{(1)}[0,1]}.
\end{equation}

Since $L_2^{(1)}[0,1]$ is a Hilbert space, we obtain
\begin{equation}
\left( {\ell,\varphi} \right) = \left\langle {\psi _\ell,\varphi}
\right\rangle \label{eq8}
\end{equation}
using the Riesz theorem for $\psi_\ell$, where $\left\langle {\psi _\ell,\varphi} \right\rangle$
is the inner product of the functions $\psi _\ell$ and $\varphi$ defined by (\ref{eq4}) and
$\varphi\in L_2^{(1)}[0,1]$, respectively.
In addition, the equality  $\|\ell\|_{L_2^{(1)*}[0,1]} = \|\psi_\ell\|_{L_2^{(1)}[0,1]}$ is achieved.
Then we obtain
\begin{equation}\label{eq9}
\left(\ell,\psi_{\ell}\right) = \|\ell\|_{L_2^{(1)*}[0,1]}^2
\end{equation}
from (\ref{eq7}). In order for  the error functional (\ref{eq3}) to be defined in the space $L_2^{(1)}[0,1]$, the condition
\begin{equation}\label{eq10}
(\ell,1)=0
\end{equation}
must be imposed which means that the quadrature formula (\ref{eq1}) is exact for any constant term.

For $\psi_{\ell}$ in (\ref{eq8}) we have
\begin{eqnarray}
&&\psi _\ell''(x)=- \bar\ell(x),         \label{eq11}\\
&&\psi _\ell'(0) =0,\ \psi _\ell'(1) =0, \label{eq12}
\end{eqnarray}
where $\bar\ell$ is the complex conjugate to $\ell$. Then the following theorem holds.

\begin{theorem}\label{Thm1}
The solution of the boundary value
problem (\ref{eq11})-(\ref{eq12}) is the extremal
function $\psi_{\ell}$ of the error functional $\ell$, expressed as
 \begin{equation}\label{eq13}
\psi_\ell(x)=-\bar\ell (x)*G_1(x)+ p_0,
 \end{equation}
where
 \begin{equation}\label{eq14}
G_1(x) = \frac{|x|}{2},
\end{equation}
$p_0=p_0^R+\i p_0^I$, a complex number, and $*$ is the convolution operation.
\end{theorem}

From Sobolev's result (see \cite{Sobolev74,SobVas}) on the extremal function of quadrature formulas in the space $L_2^{(m)}$,
we can get the statement of Theorem \ref{Thm1}, especially when $m=1$.

Next, we assume that
\begin{equation}\label{eq15}
C_{\beta}=C_{\beta}^R+\i C_{\beta}^I,
\end{equation}
where $C_{\beta}^R$ and $C_{\beta}^I$ are real numbers.
Then, using (\ref{eq10}) and (\ref{eq13}) for the norm of the error functional $\ell$ with (\ref{eq9}), we get
\begin{eqnarray*}
\|\ell\|^2&=&(\ell,\psi_{\ell})\nonumber \\
&=&\int\limits_{-\infty}^{\infty}\ell(x)\psi_{\ell}(x)\d x=-\int\limits_{-\infty}^{\infty}\ell(x)\cdot (\bar \ell(x)*G_1(x))\d x.
\end{eqnarray*}
Therefore, by direct calculation with (\ref{eq15}), we get
\begin{eqnarray}
\|\ell\|^2&=&-\Bigg[\sum\limits_{\beta=0}^N\sum\limits_{\gamma=0}^N(C_{\beta}^RC_{\gamma}^R+C_{\beta}^IC_{\gamma}^I)\ G_1(h\beta-h\gamma)\nonumber\\
&&-2\sum\limits_{\beta=0}^NC_{\beta}^R\int\limits_0^1\cos 2\pi \omega x\cdot G_1(x-h\beta)\d x\nonumber \\
&&-2\sum\limits_{\beta=0}^NC_{\beta}^I\int\limits_0^1\sin 2\pi \omega x\cdot G_1(x-h\beta)\d x\nonumber\\
&&+\int\limits_0^1\int\limits_0^1\cos[2\pi \omega(x-y)]\cdot G_1(x-y)\d x\d y\Bigg].\label{eq16}
\end{eqnarray}
Then from (\ref{eq10}) with (\ref{eq15}), we obtain the following equalities:
\begin{eqnarray}
&&\sum\limits_{\beta=0}^NC_{\beta}^R=\int\limits_0^1\cos 2\pi\omega x\ \d x,\label{eq17}\\
&&\sum\limits_{\beta=0}^NC_{\beta}^I=\int\limits_0^1\sin 2\pi\omega x\ \d x.\label{eq18}
\end{eqnarray}

Thus, we get the expression (\ref{eq16}) for the norm of the error functional (\ref{eq3}).

Further, in the next section we will solve Problem \ref{Prob1}.

\subsection{Minimization of the expression (\ref{eq16}) by coefficients $C_\beta$}

Problem \ref{Prob1} is equivalent to the problem minimizing (\ref{eq16}) in $C_\beta^R$ and $C_\beta^I$ using Lagrange method under the conditions (\ref{eq17}) and (\ref{eq18}).

Now we consider the function
\begin{eqnarray*}
&&\Psi(C_0^R,C_1^R,...,C_N^R,C_0^I,C_1^I,...,C_N^I,p_0^R,p_0^I)\\
&&=\left\| \ell\right\|^2+2p_0^R\left(\int_0^1\cos 2\pi \omega x \d x-\sum_{\beta=0}^N
C_{\beta}^R\right)+2p_0^I\left(\int_0^1\sin 2\pi \omega x \d x-\sum_{\beta=0}^N
C_{\beta}^I\right).
\end{eqnarray*}
By making the partial derivatives of
$\Psi$ with respect to $C_\beta^R$, $C_{\beta}^I$, $(\beta =\overline{0,N})$,
$p_0^R$ and $p_0^I$ equal to zero, we get the following
system of linear equations:
\begin{eqnarray}
&&\sum\limits_{\gamma=0}^N C_\gamma^R G_1(h\beta   - h\gamma) +
p_0^R= \int\limits_0^1\cos 2\pi \omega x G_1(x-h\beta)\d x, \beta=0,...,N,\label{eq19}\\
&&\sum\limits_{\gamma=0}^N C_\gamma^R
=\int\limits_0^1\cos 2\pi \omega x\ \d x,\label{eq20}\\
&&\sum\limits_{\gamma=0}^N C_\gamma^I G_1(h\beta   - h\gamma) +
p_0^I = \int\limits_0^1\sin 2\pi \omega x G_1(x-h\beta)\d x, \beta=0,...,N,\label{eq21} \\
&&\sum\limits_{\gamma=0}^N C_\gamma^I
=\int\limits_0^1\sin 2\pi \omega x\ \d x.\label{eq22}
\end{eqnarray}
We multiple  both sides of (\ref{eq21}) and (\ref{eq22}) by $\i$ and add these to
(\ref{eq19}) and (\ref{eq20}), respectively, to obtain a system
of $(N+2)$ linear equations with $(N+2)$ unknowns $C_{\gamma}$, $\gamma=0,1,...,N$, and $p_0$:
\begin{eqnarray}
&&\sum\limits_{\gamma=0}^N C_\gamma G_1(h\beta   - h\gamma) +
p_0 = \int_0^1\E^{2\pi \i\omega x} G_1(x - h\beta)\d x,\ \beta=0,...,N,\label{eq23} \\
&&\sum\limits_{\gamma=0}^N C_\gamma=\int\limits_0^1\E^{2\pi \i \omega x}\d x,\label{eq24}
\end{eqnarray}
where $G_1(x)$ is defined in (\ref{eq14}). The system (\ref{eq23})-(\ref{eq24}) has a unique solution. The uniqueness of the solution
of this system can be proved by the uniqueness of the solution of the system (3.1)-(3.2) in \cite{ShadHayAkhm15}.
The solution of the system (\ref{eq23})-(\ref{eq24}) provides the minimum of $\left\|\ell\right\|^2$
at $C_{\beta}=\mathring{C}_{\beta}$. The quadrature formula of the form (\ref{eq1}) with coefficients
$\mathring{C}_{\beta}$ is called \emph{the optimal quadrature formula} in the sense of Sard, and
$\mathring{C}_{\beta}$ are said to be \emph{the optimal coefficients}.
For convenience, the optimal coefficients $\mathring{C}_{\beta}$ will be denoted as $C_{\beta}$.

The purpose of this section is to obtain an analytic solution for the system (\ref{eq23})-(\ref{eq24}).
To do this, we  use the concept of discrete argument functions and
operations. The theory of discrete argument functions is
given in \cite{Sobolev74,SobVas}. We give
the definition for the function of discrete argument.
Suppose that nodes $x_\beta$ has uniform spacing (i.e.,
$x_\beta=h\beta,$ $h$ is a small positive parameter), and functions $\varphi(x)$ and
$\psi(x)$ are complex-valued and defined on the real line $\mathbb{R}$ or on an interval of $\mathbb{R}$.

The function $\varphi (h\beta )$ is {\it a
function of discrete argument } {if it is given on some set
of integer values of} $\beta$.
{\it The inner product} of two discrete argument
functions $\varphi(h\beta )$ and $\psi (h\beta )$ is given by
$$
\left[ {\varphi(h\beta),\psi(h\beta) } \right] =
\sum\limits_{\beta  =  - \infty }^\infty  {\varphi (h\beta ) \cdot
\bar \psi (h\beta )},
$$
{if the series on the right hand side of the last equality
converges absolutely.}
{\it The convolution} {of two
functions $\varphi(h\beta )$ and $\psi (h\beta )$ is the inner
product}
$$
\varphi (h\beta )*\psi (h\beta ) = \left[ {\varphi (h\gamma ),\psi
(h\beta  - h\gamma )} \right] = \sum\limits_{\gamma  =  - \infty
}^\infty  {\varphi (h\gamma ) \cdot \bar\psi (h\beta  - h\gamma )}.
$$

We also use the discrete analogue $D_1(h\beta)$ for the operator $\d^2/\d x^2$, that satisfies
\begin{equation}\label{eq25}
hD_1(h\beta)*G_1(h\beta)=\delta_{\d}(h\beta),
\end{equation}
where $G_1(h\beta)=\frac{|h\beta|}{2}$,
$\delta_{\d}(h\beta)$ is equal to 0 when $\beta\neq 0$, and 1 when $\beta=0$.

It should be noted that the discrete analogue $D_m(h\beta)$ of the differential operator $\d^{2m}/\d x^{2m}$ was first
introduced and investigated by Sobolev \cite{Sobolev74,SobVas} and it was constructed in \cite{Shad85}. In particular,
from the results of \cite{Shad85} for $m=1$, the following are obtained.

\begin{theorem}\label{Thm2}
The discrete analogue $D_1(h\beta)$ to the operator $\d^2/\d x^2$ satisfying (\ref{eq25})
has the form
\begin{equation}\label{eq26}
D_1(h\beta)=\frac{1}{h^2}\left\{
\begin{array}{rl}
0,&\ |\beta|\geq 2,\\
1,&\ |\beta|=1,\\
-2,& \  \beta=0
\end{array}
\right.
\end{equation}
and satisfies
\begin{equation}\label{eq27}
D_1(h\beta)*1=0,\ \ D_1(h\beta)*(h\beta)=0.
\end{equation}
\end{theorem}

Now we return to our problem.

We regard the coefficients $C_{\beta}$ as a discrete argument function and assume $C_{\beta}=0$ for
$\beta=-1,-2,...$ and $\beta=N+1,N+2,...$. Then, considering the above definitions,
we rewrite the system (\ref{eq23})-(\ref{eq24}) in the convolution form as
\begin{eqnarray}
&&C_{\beta}*G_1(h\beta) +
p_0 = f_1(h\beta),\ \ \beta  =0,1,...,N, \label{eq28} \\
&&\sum\limits_{\beta=0}^N C_{\beta}=g_0,\label{eq29}
\end{eqnarray}
where
\begin{eqnarray}
f_1(h\beta)&=&-\frac{h\beta}{2(2\pi\i \omega)}(\E^{2\pi\i \omega}+1)\nonumber \\
&&+\frac{1}{2(2\pi \i \omega)^2}\left(2\E^{2\pi\i \omega h\beta}+(2\pi\i \omega-1)\E^{2\pi\i \omega}-1\right),\label{eq30} \\
g_0 &=&\frac{1}{2\pi\i \omega}(\E^{2\pi\i \omega}-1)\label{eq31}
\end{eqnarray}
and $G_1(x)$ is defined by (\ref{eq14}).

Now we have the following problem.

\begin{problem}\label{Prob2}
Find $C_{\beta}$, $\beta=0,1,..,N$, and $p_0$ satisfying the system (\ref{eq28})-(\ref{eq29}) for given $f_1(h\beta)$ and $g_0$.
\end{problem}

Note that Problem \ref{Prob2} is equivalent to Problem \ref{Prob1}.
The main result of this section is as follows.

\begin{theorem}\label{Thm3}
For $\omega\in \mathbb{R}$ with $\omega\neq 0$, coefficients of the optimal quadrature formulas of the form (\ref{eq1}) in the sense of Sard  in the space $L_2^{(1)}[0,1]$ have the form
\begin{equation}\label{eq32}
\begin{array}{rcl}
{C}_0&=&\displaystyle h\ \frac{(1+2\pi\i \omega h-\E^{2\pi\i \omega h})}{(2\pi \omega h)^2},\\
{C}_\beta&=&\displaystyle h\ \frac{2(1-\cos 2\pi\omega h)}{(2\pi \omega h)^2}\ \E^{2\pi\i \omega h\beta},\ \beta=1,2,...,N-1,\\
{C}_N&=&\displaystyle h\ \frac{(1-2\pi\i \omega h-\E^{-2\pi\i \omega h})}{(2\pi \omega h)^2}\ \E^{2\pi \i \omega}.
\end{array}
\end{equation}
In addition, for the square of the norm of the error functional (\ref{eq3}) of the optimal quadrature formula (\ref{eq1}) in the space
$L_2^{(1)*}[0,1]$, the following holds:
\begin{equation}\label{eq33}
\left\|\mathring{\ell}\right\|_{L_2^{(1)*}}^2=\frac{1}{(2\pi\omega)^2}\left(1-\frac{2(1-\cos 2\pi\omega h)}{(2\pi \omega h)^2}\right).
\end{equation}
\end{theorem}

\emph{Proof.} We consider a discrete argument function
\begin{equation}\label{eq34}
u_1(h\beta)=C_\beta*G_1(h\beta)+p_0.
\end{equation}
Then, considering (\ref{eq25}) and (\ref{eq27}), we have
\begin{equation}\label{eq35}
C_{\beta}=hD_1(h\beta)*u_1(h\beta).
\end{equation}
Calculating the convolution (\ref{eq35}) requires the representation of the function $u_1(h\beta)$ for all integer values of $\beta$. From (\ref{eq28}) we have
\begin{equation}\label{eq36}
u_1(h\beta)=f_1(h\beta)\mbox{ for }\beta=0,1,...,N.
\end{equation}
Now we need to find the representation of $u_1(h\beta)$ for $\beta<0$ and $\beta>N$.
Using (\ref{eq14}) and (\ref{eq29}) for $\beta\leq 0$ and $\beta\geq N$, respectively, we get
\begin{equation}\label{eq37}
u_1(h\beta)=
\left\{
\begin{array}{ll}
-\frac{h\beta}{2}\ g_0+\frac{1}{2}\sum\limits_{\gamma=0}^NC_{\gamma}h\gamma+p_0,&\ \beta\leq 0,\\
\frac{h\beta}{2}\ g_0-\frac{1}{2}\sum\limits_{\gamma=0}^NC_{\gamma}h\gamma+p_0,& \ \beta\geq N,
\end{array}
\right.
\end{equation}
where $g_0$ is defined as (\ref{eq31}), and $\sum\limits_{\gamma=0}^NC_{\gamma}h\gamma$ and $p_0$ are unknowns.
Then from the last two equalities when $\beta=0$ and $\beta=N$, we get the following system of two linear equations
for these unknowns:
\begin{eqnarray*}
&&p_0+\frac{1}{2}\sum\limits_{\gamma=0}^NC_{\gamma}h\gamma=f_1(0),\\
&&p_0-\frac{1}{2}\sum\limits_{\gamma=0}^NC_{\gamma}h\gamma+\frac{1}{2}g_0=f_1(1).
\end{eqnarray*}
Therefore, solving this system using (\ref{eq30}) and (\ref{eq31}), we get
\begin{eqnarray}
p_0&=&0,
\label{eq38}\\
\sum\limits_{\gamma=0}^NC_{\gamma}h\gamma&=&\frac{\E^{2\pi \i \omega}}{2\pi \i \omega}-\frac{\E^{2\pi \i \omega}-1}{(2\pi \i \omega)^2}.
\label{eq39}
\end{eqnarray}
With (\ref{eq38}) and (\ref{eq39}) in mind, the combination of (\ref{eq36}) and (\ref{eq37}) results in
\begin{equation*}
u_1(h\beta)=
\left\{
\begin{array}{ll}
-\frac{h\beta}{2}g_0+\frac{(2\pi \i \omega-1)\E^{2\pi \i \omega}+1}{2(2\pi\i \omega)^2},&\ \beta\leq 0,\\
f_1(h\beta),&\ 0\leq \beta\leq N,\\
\frac{h\beta}{2} g_0-\frac{(2\pi \i \omega-1)\E^{2\pi \i \omega}+1}{2(2\pi\i \omega)^2},& \ \beta\geq N.
\end{array}
\right.
\end{equation*}
The analytic formulas (\ref{eq32}) is now obtained from (\ref{eq35}) by taking into account (\ref{eq26}) and (\ref{eq27}),
using the last representation of $u_1(h\beta)$,
and  by direct calculation of the optimal coefficients ${C}_\beta,$ $\beta=0,1,...,N$.

Now we are going to get (\ref{eq33}). We rewrite (\ref{eq16}) in the
following form:
\begin{eqnarray}
\|\mathring{\ell}\|^2&=&-\Bigg[\sum\limits_{\beta=0}^NC_\beta^R\left(\sum\limits_{\gamma=0}^NC_\gamma^R
G_1(h\beta-h\gamma)- \int\limits_0^1\cos 2\pi \omega x\  G_1(x-h\beta)\d x\right)\nonumber\\
&&+\sum\limits_{\beta=0}^NC_\beta^I\left(\sum\limits_{\gamma=0}^NC_\gamma^I
G_1(h\beta-h\gamma)- \int\limits_0^1\sin 2\pi \omega x\  G_1(x-h\beta)\d x\right)\nonumber\\
&& -\sum\limits_{\beta=0}^NC_\beta^R \int\limits_0^1\cos 2\pi \omega x\  G_1(x-h\beta)\d x
-\sum\limits_{\beta=0}^NC_\beta^I \int\limits_0^1\sin 2\pi \omega x\  G_1(x-h\beta) \d x\nonumber\\
&&+\int\limits_0^1\int\limits_0^1\cos[2\pi \omega(x-y)]G_1(x-y)\d x\d y\Bigg].
\label{eq40}
\end{eqnarray}
Since $p_0=p_0^R+\i p_0^I$, considering (\ref{eq38}), we have
$$
p_0^R=0\mbox{ and } p_0^I=0.
$$
Therefore, these two last equalities are used in (\ref{eq19}) and (\ref{eq21}) to obtain
$$
\sum\limits_{\gamma=0}^NC_\gamma^R
G_1(h\beta-h\gamma)- \int\limits_0^1\cos 2\pi \omega x\  G_1(x-h\beta)\d x=0,\ \beta=0,...,N
$$
and
$$
\sum\limits_{\gamma=0}^NC_\gamma^I
G_1(h\beta-h\gamma)- \int\limits_0^1\sin 2\pi \omega x\  G_1(x-h\beta)\d x=0,\ \beta=0,...,N.
$$
Then the expression (\ref{eq40}) for $\|\mathring{\ell}\|^2$ takes the form
\begin{eqnarray*}
\|\mathring{\ell}\|^2&=&\sum\limits_{\beta=0}^NC_\beta^R \int\limits_0^1\cos 2\pi \omega x\  G_1(x-h\beta)\d x
+\sum\limits_{\beta=0}^NC_\beta^I \int\limits_0^1\sin 2\pi \omega x\  G_1(x-h\beta) \d x\nonumber\\
&&-\int\limits_0^1\int\limits_0^1\cos[2\pi \omega(x-y)]G_1(x-y)\d x\d y.
\end{eqnarray*}
Therefore calculating the definite integrals, keeping (\ref{eq15}) in mind and using (\ref{eq32}),
we get (\ref{eq33}) after some simplifications. Theorem \ref{Thm3} has been proved.~$\Box$

We note that in Theorem \ref{Thm3}, the formulas for the optimal coefficients ${C}_\beta$ are decomposed into two parts:
real and imaginary parts. Therefore from the formulas (\ref{eq32}) of Theorem \ref{Thm3}, we get the following results.

\begin{corollary}\label{Cor1}
{For $\omega\in \mathbb{R}$ with $\omega\neq 0$,  coefficients of the optimal quadrature formula of the form
$$
\int\limits_0^1\cos 2\pi \omega x\cdot \varphi(x)\d x\cong \sum\limits_{\beta=0}^NC_{\beta}^R\varphi(h\beta)
$$
in the sense of Sard in $L_2^{(1)}[0,1]$  have the form}
\begin{eqnarray*}
{C}_0^R&=&h\ \frac{1-\cos 2\pi \omega h}{(2\pi\omega h)^2},\\
{C}_\beta^R&=&h\ \frac{2(1-\cos 2\pi \omega h)}{(2\pi \omega h)^2}\ \cos 2\pi \omega h\beta,\ \ \beta=1,2,...,N-1,\\
{C}_N^R&=&h\ \frac{(1-\cos 2\pi \omega h)\cos 2\pi \omega +(2\pi\omega h-\sin 2\pi \omega h)\sin 2\pi \omega}{(2\pi\omega h)^2}.\\
\end{eqnarray*}
\end{corollary}

\begin{corollary}\label{Cor2}
{For $\omega\in \mathbb{R}$ with $\omega\neq 0$, coefficients of the optimal quadrature formula of the form
$$
\int\limits_0^1\sin 2\pi \omega x\cdot \varphi(x)\d x\cong \sum\limits_{\beta=0}^NC_{\beta}^I\varphi(h\beta)
$$
in the sense of Sard
in $L_2^{(1)}[0,1]$   have the form}
\begin{eqnarray*}
{C}_0^I&=&h\ \frac{2\pi \omega h-\sin 2\pi \omega h}{(2\pi \omega h)^2},\\
{C}_\beta^I&=&h\ \frac{2(1-\cos 2\pi \omega h)}{(2\pi \omega h)^2}\ \sin{2\pi \omega h\beta},\ \ \beta=1,2,...,N-1,\\
{C}_N^I&=&h\ \frac{(1-\cos 2\pi \omega h)\sin 2\pi \omega -(2\pi\omega h-\sin 2\pi \omega h)\cos 2\pi \omega}{(2\pi\omega h)^2}.\\
\end{eqnarray*}
\end{corollary}

It is easy to see that for $\omega\to 0$ Sard's following result \cite{Sard} on the optimality of the trapezoidal
quadrature formula in $L_2^{(1)}[0,1]$ is obtained from Theorem~\ref{Thm3}.

\begin{corollary}\label{Cor3}
Coefficients of the optimal quadrature formula of the form
\begin{equation}\label{eq41}
\int\limits_0^1 \varphi(x)\d x\cong \sum\limits_{\beta=0}^NC_{\beta}\varphi(h\beta)
\end{equation}
in the space $L_2^{(1)}[0,1]$ have the form
\begin{equation*}
\begin{array}{rcl}
{C}_0&=&\displaystyle \frac{h}{2},\\
{C}_\beta&=&\displaystyle h,\ \beta=1,2,...,N-1,\\
{C}_N&=&\displaystyle \frac{h}{2}
\end{array}
\end{equation*}
and for the norm of the error functional of the optimal quadrature formula (\ref{eq41}) in the space
$L_2^{(1)*}[0,1]$, the following holds
\begin{equation*}
\left\|\mathring{\ell}\right\|_{L_2^{(1)*}}^2=\frac{h^2}{12}.
\end{equation*}
\end{corollary}

In addition,  for $\omega h\in \mathbb{Z}$ with $\omega\neq 0$ we obtain
the following corollary from (\ref{eq32}) and (\ref{eq33}).

\begin{corollary}\label{Cor4} For $\omega h\in \mathbb{Z}$ with $\omega\neq 0$,
coefficients of the optimal quadrature formula of the form (\ref{eq1})
in the sense of Sard in the space $L_2^{(1)}[0,1]$ have the form
\begin{equation*}
\begin{array}{rcl}
{C}_0&=&\displaystyle \frac{1}{2\pi\i \omega},\\
{C}_\beta&=&0,\ \beta=1,2,...,N-1,\\
{C}_N&=&\displaystyle -\frac{1}{2\pi\i \omega}
\end{array}
\end{equation*}
and for the norm of the error functional (\ref{eq3}) of the optimal quadrature formula (\ref{eq1}) in the space
$L_2^{(1)*}[0,1]$, the following holds:
$$
\|\mathring{\ell}\|_{L_2^{(1)*}}^2=\frac{1}{(2\pi \omega)^2},
$$
i.e., the convergence order of the optimal quadrature formula of the form (\ref{eq1}) is $O(|\omega|^{-1})$ for $\omega h\in \mathbb{Z}$ with $\omega\neq 0$.
\end{corollary}

{\it Remark 1} It should be noted that for a fixed $\omega$, we obtain
$$
\|\mathring{\ell}\|^2=\frac{1}{12}h^2-\frac{1}{90}\pi^2 \omega^2 h^4+\frac{1}{1260}\pi^4 \omega^4 h^6+O(h^8),
$$
from (\ref{eq33}), i.e., the convergence order of the optimal quadrature formula of the form (\ref{eq1}) is $O(h)$.

{\it Remark 2} In particular,  in the case $\omega\in \mathbb{Z}$ with $\omega\neq 0$,
the results of \cite{BolHayShad16} and of Section 6 of \cite{BolHayShad17} are obtained from Theorem \ref{Thm3}.

{\it Remark 3} The equality (\ref{eq39}) means that the optimal quadrature formula of the form (\ref{eq1}) with coefficients (\ref{eq32})
is exact to $\varphi(x)=x$ because
$$
\int\limits_0^1\E^{2\pi\i \omega x}x\ \d x=\frac{\E^{2\pi \i \omega}}{2\pi \i \omega}-\frac{\E^{2\pi \i \omega}-1}{(2\pi \i \omega)^2}.
$$
The equality (\ref{eq39}) together with (\ref{eq29}) provides the exactness of our optimal quadrature formula for all
linear functions. Therefore, for functions with a continuous second derivative, the convergence order
of the optimal quadrature formula (\ref{eq1}) with coefficients (\ref{eq32})
is concluded as $O(h^2)$.

\section{Optimal quadrature formula for the interval [a,b]}
\setcounter{equation}{0}

Here, optimal quadrature formulas for the interval $[a,b]$ are obtained
by a linear transform from the results of the previous section.

We consider the construction of the optimal quadrature formula of the form
\begin{equation}\label{eq3.1}
\int\limits_a^b\E^{2\pi \i \omega x}\varphi(x)\ \d x\cong \sum\limits_{\beta=0}^NC_{\beta,\omega}[a,b]\varphi(x_\beta)
\end{equation}
in the Sobolev space
$L_2^{(1)}[a,b]$. Here $C_{\beta,\omega}[a,b]$ are coefficients,
$x_\beta=h\beta+a$ $(\in [a,b])$ are the nodes of the formula (\ref{eq3.1}),
$\omega\in \mathbb{R}$, $\i^2=-1$, and $h=\frac{b-a}{N}$ for $N\in \mathbb{N}$.

Now, by a linear transformation $x=(b-a)y+a$, where $0\leq y\leq 1$, we obtain
$$
\int\limits_a^b\E^{2\pi \i \omega x}\varphi(x)\ \d x=(b-a)\E^{2\pi \i \omega a}\int\limits_0^1\E^{2\pi \i \omega (b-a)y}\varphi((b-a)y+a)\d y.
$$
Finally, by applying Theorem \ref{Thm3} and Corollary \ref{Cor3} to the integral on the right-hand side of the last equality, we have the following
main result of the present work.

\begin{theorem}\label{Thm4}
For $\omega\in \mathbb{R}$ with $\omega\neq 0$, coefficients of the optimal quadrature formula of the form
\begin{equation}\label{eq3.2}
\int\limits_a^b\E^{2\pi \i \omega x}\varphi(x)\ \d x\cong \sum\limits_{\beta=0}^NC_{\beta,\omega}[a,b] \varphi(h\beta+a)
\end{equation}
in the sense of Sard in the space $L_2^{(1)}[a,b]$ have the form
\begin{equation}\label{eq3.3}
\begin{array}{rcl}
{C}_{0,\omega}[a,b]&=&\displaystyle h\ \frac{(1+2\pi\i \omega h-\E^{2\pi\i \omega h})}{(2\pi \omega h)^2}\ \E^{2\pi \i \omega a},\\
{C}_{\beta,\omega}[a,b]&=&\displaystyle h\ \frac{2(1-\cos 2\pi\omega h)}{(2\pi \omega h)^2}\ \E^{2\pi\i \omega (h\beta+a)},\ \beta=1,2,...,N-1,\\
{C}_{N,\omega}[a,b]&=&\displaystyle h\ \frac{(1-2\pi\i \omega h-\E^{-2\pi\i \omega h})}{(2\pi \omega h)^2}\ \E^{2\pi \i \omega b},
\end{array}
\end{equation}
and for $\omega=0$, the coefficients take the form
\begin{equation}\label{eq3.4}
\begin{array}{rcl}
{C}_{0,0}[a,b]&=&\displaystyle \frac{h}{2},\\
{C}_{\beta,0}[a,b]&=&\displaystyle h,\ \beta=1,2,...,N-1,\\
{C}_{N,0}[a,b]&=&\displaystyle \frac{h}{2},
\end{array}
\end{equation}
where $h=\frac{b-a}{N}$.
\end{theorem}

The monomials $x^{\alpha}$ for $\alpha=0,1,2,...$ are now considered as a function $\varphi$ in the integral in the left-hand side of (\ref{eq3.2}).
Then we get
\begin{eqnarray}
g_{\alpha,\omega}[a,b]&=&\int\limits_a^b\E^{2\pi \i\omega x}x^{\alpha}\d x\nonumber\\
&=&
\left\{
\begin{array}{ll}
\sum\limits_{k=0}^{\alpha-1}\frac{(-1)^k\alpha!}{(\alpha-k)!\ (2\pi \i \omega)^{k+1}}
\left(\E^{2\pi \i\omega b}b^{\alpha-k}-\E^{2\pi \i\omega a}a^{\alpha-k} \right)& \\
\qquad\qquad\qquad+\frac{(-1)^{\alpha}\alpha!}{(2\pi \i \omega)^{\alpha+1}}\left(\E^{2\pi \i\omega b}-\E^{2\pi \i\omega a} \right)&\mbox{ for } \omega\neq 0,\\
\frac{1}{\alpha+1}(b^{\alpha+1}-a^{\alpha+1})&\mbox{ for }\omega=0,
\end{array}
\right.\label{eq3.5}
\end{eqnarray}
where $\alpha=0,1,2,...$.

\emph{Remark 4} We find from Remark 3 that the optimal quadrature formula (\ref{eq3.2}) is exact for all linear functions, that is
the following equalities hold:
\begin{eqnarray}
g_{0,\omega}[a,b]&=&\sum\limits_{\beta=0}^N{C}_{\beta,\omega}[a,b],\label{eq3.6}\\
g_{1,\omega}[a,b]&=&\sum\limits_{\beta=0}^N{C}_{\beta,\omega}[a,b]\cdot(h\beta+a),\label{eq3.7}
\end{eqnarray}
where $g_{0,\omega}[a,b]$ and $g_{1,\omega}[a,b]$ are defined by (\ref{eq3.5}) when $\alpha=0$ and $\alpha=1$, respectively.

\section{Approximation of Fourier transforms by optimal quadrature formula}
\setcounter{equation}{0}

Here we consider some numerical results confirming the theoretical results of the previous sections.
The present section consists of two parts.
In the first part, using  the optimal quadrature formula (\ref{eq3.2}),
we approximate the integrals
$$
g_{\alpha,\omega}[-1,1]=\int\limits_{-1}^1\E^{2\pi \i \omega x}x^{\alpha}\d x, \ \ \alpha=0,1,2,
$$
where
\begin{eqnarray}
g_{0,\omega}[-1,1]&=&\left\{
\begin{array}{ll}
\frac{1}{\pi \omega}\sin 2\pi \omega,&\ \omega\neq 0,\\
2,&\ \omega=0,
\end{array}
\right. \label{eq4.1}\\
g_{1,\omega}[-1,1]&=&\left\{
\begin{array}{ll}
\frac{2\i}{(2\pi \omega)^2}(\sin 2\pi \omega-2\pi \omega\cos 2\pi\omega),&\ \omega\neq 0,\\
0,&\ \omega=0,
\end{array}
\right. \label{eq4.2}\\
g_{2,\omega}[-1,1]&=&\left\{
\begin{array}{ll}
\frac{4}{(2\pi \omega)^3}((2\pi^2 \omega^2-1)\sin 2\pi \omega+2\pi \omega\cos 2\pi\omega),&\ \omega\neq 0,\\
\frac{2}{3},&\ \omega=0.
\end{array}
\right. \label{eq4.3}
\end{eqnarray}

In the second part, using the given function $f$, the optimal quadrature formula (\ref{eq3.2}) is applied to the approximation of the  Fourier Transforms
\begin{eqnarray}
F(\omega)&=&\int\limits_{-\infty}^{\infty}\E^{-2\pi \i\omega x}f(x)\d x,\label{eq4.4} \\
f(x)&=&\int\limits_{-\infty}^{\infty}\E^{2\pi \i\omega x}F(\omega)\d \omega \label{eq4.5}
\end{eqnarray}
and thereby resulting in the approximate reconstruction of the function $f$.

\subsection{Approximation of the Fourier integral }

The error of the optimal quadrature formula (\ref{eq3.2})  is denoted by
$$
R_{\varphi,\omega}[a,b]=\int\limits_a^b\E^{2\pi \i\omega x}\varphi(x)\d x-
\sum\limits_{\beta=0}^N{C}_{\beta,\omega}[a,b]\varphi(h\beta+a).
$$

Consider the functions $f_{\alpha}$, $\alpha=0,1,2$, obtained
by extending monomials $x^{\alpha}$, $\alpha=0,1,2$, respectively, with zeros outside the interval $[-1,1]$, i.e.,
we have
\begin{equation}\label{eq4.7}
f_{\alpha}(x)=
\left\{
\begin{array}{ll}
x^{\alpha}& \mbox{ for } \ x\in [-1,1],\\
0& \mbox{ otherwise}
\end{array}
\right.\ \mbox{ for } \alpha=0,1,2.
\end{equation}

First the optimal quadrature formula (\ref{eq3.2}) is applied to approximate the integrals
$$
\int\limits_{-1}^1\E^{2\pi \i \omega x}f_{\alpha}(x)\d x
$$
for functions $f_{\alpha}$, $\alpha=0,1,2$, in (\ref{eq4.7}).

\begin{figure}[h]
   \includegraphics[width=\textwidth]{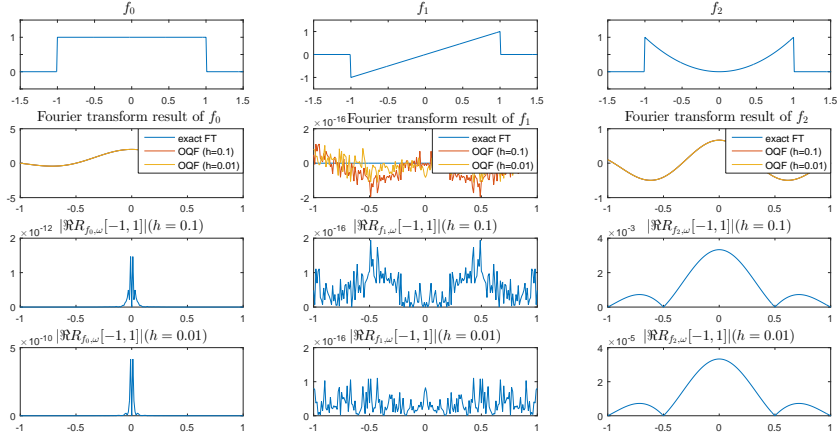}
\caption{Graphs of functions $f_{\alpha}$, $\alpha=0,1,2$, defined by (\ref{eq4.7}) (the first row), graphs of the exact Fourier transforms for functions
$f_{\alpha}$, $\alpha=0,1,2$ and their approximations by OQF applied for the interval $[-1,1]$ with steps $h=0.1$ and $h=0.01$ (the second row)
as well as graphs of $|\Re{R}_{f_{\alpha},\omega}[-1,1]|$, $\alpha=0,1,2$ when $h=0.1$ (the third row) and $h=0.01$ (the forth row) with $\omega\in [-1,1]$.}
\label{Fig4.1}
\end{figure}

With (\ref{eq3.6}) and (\ref{eq3.7}), for functions $f_0$ and $f_1$, we obtain
$$
g_{0,\omega}[-1,1]=\sum\limits_{\beta=0}^N{C}_{\beta,\omega}[-1,1],\ \
g_{1,\omega}[-1,1]=\sum\limits_{\beta=0}^N{C}_{\beta,\omega}[-1,1]\cdot(h\beta-1),
$$
where $h=2/N$. Therefore, $R_{f_0,\omega}[-1,1]=0$ and $R_{f_1,\omega}[-1,1]=0$. The graphs for the absolute values of the real parts of the errors
$R_{f_0,\omega}[-1,1]$ and $R_{f_1,\omega}[-1,1]$ in the first and the second columns of Figure \ref{Fig4.1}
confirm the real parts of the last equalities numerically. For the function $f_2$, there is the error
\begin{equation}\label{eq4.6}
R_{f_2,\omega}[-1,1]=g_{2,\omega}[-1,1]-\sum\limits_{\beta=0}^N{C}_{\beta,\omega}[-1,1]\ (h\beta-1)^2,
\end{equation}
where ${C}_{\beta,\omega}[-1,1]$ are defined by (\ref{eq3.3}) and (\ref{eq3.4}) with $h=2/N$.
For $h=0.1$ and $h=0.01$ with $\omega\in [-1,1]$ the graphs for absolute values of the real part of the error (\ref{eq4.6}) are shown in the third column of
Figure \ref{Fig4.1}. From the graphs in the third and forth rows of the third column of Figure \ref{Fig4.1}, we can see that the error (\ref{eq4.6}) of the optimal quadrature formula (\ref{eq3.2}) for the function $f_2(x)=x^2$ is $O(h^2)$. This statement confirms Remark 3 numerically.

\begin{figure}[h]
   \includegraphics[width=\textwidth]{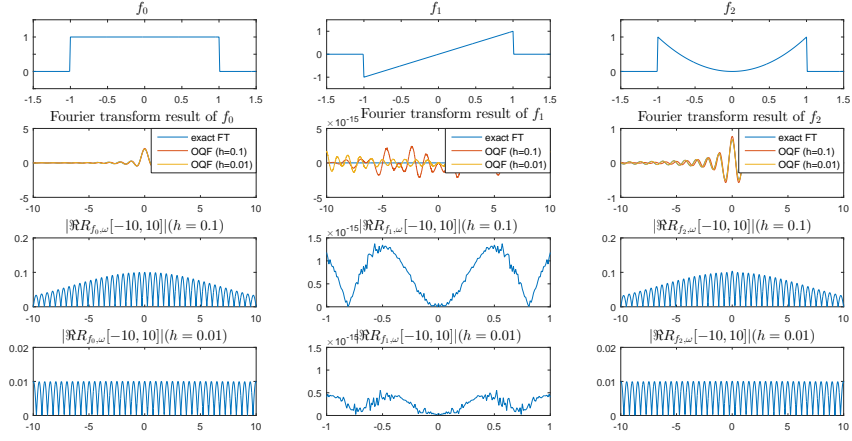}
\caption{Graphs of functions $f_{\alpha}$, $\alpha=0,1,2$, defined by (\ref{eq4.7}) (the first row), graphs of the exact Fourier transforms for functions
$f_{\alpha}$, $\alpha=0,1,2$ and their approximations by OQF applied for the interval $[-10,10]$ with steps $h=0.1$ and $h=0.01$ (the second row)
as well as graphs of $|\Re{R}_{f_{\alpha},\omega}[-10,10]|$, $\alpha=0,1,2$ when $h=0.1$ (the third row) and $h=0.01$ (the forth row) with $\omega\in [-10,10]$.}
\label{Fig4.2}
\end{figure}

\begin{figure}[h]
   \includegraphics[width=\textwidth]{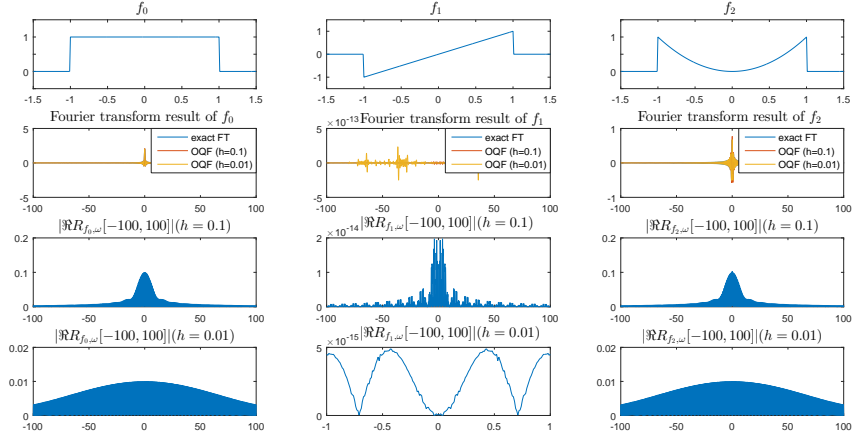}
\caption{Graphs of functions $f_{\alpha}$, $\alpha=0,1,2$, defined by (\ref{eq4.7}) (the first row), graphs of the exact Fourier transforms for functions
$f_{\alpha}$, $\alpha=0,1,2$ and their approximations by OQF applied for the interval $[-100,100]$ with steps $h=0.1$ and $h=0.01$ (the second row)
as well as graphs of $|\Re{R}_{f_{\alpha},\omega}[-100,100]|$, $\alpha=0,1,2$ when $h=0.1$ (the third row) and $h=0.01$ (the forth row) with $\omega\in [-100,100]$.}
\label{Fig4.3}
\end{figure}

Then for an interval $[a,b]$ containing the interval $[-1,1]$, the error $R_{f_{\alpha},\omega}[a,b]$, $\alpha=0,1,2$,
of the optimal quadrature formula (\ref{eq3.2}) corresponding to the functions (\ref{eq4.7}) takes the form
\begin{eqnarray}
R_{f_{\alpha},\omega}[a,b]&=&\int\limits_a^b\E^{2\pi \i\omega x}f_{\alpha}(x)\d x-
\sum\limits_{\beta=0}^N{C}_{\beta,\omega}[a,b]f_{\alpha}(h\beta+a)\nonumber\\
&=&\int\limits_{-1}^1\E^{2\pi \i\omega x}x^{\alpha}\d x-
\sum\limits_{\beta=0}^N{C}_{\beta,\omega}[a,b]f_{\alpha}(h\beta+a)\nonumber\\
&=&g_{\alpha,\omega}[-1,1]- \sum\limits_{\beta=0}^N{C}_{\beta,\omega}[a,b]f_{\alpha}(h\beta+a),\label{eq4.8}
\end{eqnarray}
where $g_{\alpha,\omega}[-1,1]$, $\alpha=0,1,2$, are defined by (\ref{eq4.1})-(\ref{eq4.3}).
We provide numerical results for intervals $[-10,10]$ and $[-100,100]$.
In these intervals, from (\ref{eq4.8}) for the errors of the optimal quadrature formula (\ref{eq3.2}), we get
\begin{eqnarray*}
R_{f_{\alpha},\omega}[-10,10]=
g_{\alpha,\omega}[-1,1]-\sum\limits_{\beta=0}^N{C}_{\beta,\omega}[-10,10]f_{\alpha}(h\beta-10)
\end{eqnarray*}
and
\begin{eqnarray*}
R_{f_{\alpha},\omega}[-100,100]=
g_{\alpha,\omega}[-1,1]-\sum\limits_{\beta=0}^N{C}_{\beta,\omega}[-100,100]f_{\alpha}(h\beta-100),
\end{eqnarray*}
respectively. For functions $f_{\alpha}$, $\alpha=0,1,2$, defined by (\ref{eq4.7}),
Figures \ref{Fig4.2} and \ref{Fig4.3} show the graphs of $f_{\alpha}$
(the first rows), graphs of the exact Fourier transforms for functions
$f_{\alpha}$ and their approximations by OQF applied for the intervals $[-10,10]$ and $[-100,100]$ with steps $h=0.1$ and $h=0.01$
(the second rows) as well as graphs of $|\Re{R}_{f_{\alpha},\omega}[-10,10]|$ and $|\Re{R}_{f_{\alpha},\omega}[-100,100]|$ when $h=0.1$ (the third rows) and $h=0.01$ (the forth rows).

We note that the functions $f_{\alpha}$ defined by (\ref{eq4.7}) are piecewise continuous and
do not belong to the space $L_2^{(1)}[a,b]$ when the interval $[a,b]$ contains the interval $[-1,1]$ and wider than it.
Nevertheless, from the numerical results in the first and the third columns of Figures \ref{Fig4.2} and \ref{Fig4.3},
we conclude that the convergence order of the optimal quadrature formula (\ref{eq3.2}) for these functions is\linebreak $O((h^{-1}+|\omega|)^{-1})$.

Note that the real part of the function $g_{1,\omega}[-1,1]$ is zero. Due to the symmetry of the considered intervals $[-1,1]$, $[-10,10]$ and $[-100,100]$ and the oddness of the function $f_1$ as well as the evenness of the optimal coefficients, the real part of the corresponding quadrature sum is also zero.
This means that the absolute values of the real part of the error $R_{f_1,\omega}[a,b]$ are zero.  This assertion confirms the numerical results (machine zero) given in the second columns of Figures \ref{Fig4.1}, \ref{Fig4.2} and \ref{Fig4.3}.

It is easy to see if the error of the optimal quadrature formula (\ref{eq3.2}) is less than the error of the Discrete Fourier Transform for the
integral $\int_{a}^b\E^{2\pi \i \omega x}\varphi(x)\d x$.

\subsection{Reconstruction of a function using approximate direct and inverse Fourier transforms}

It is known that when complete continuous X-ray data are available then CT image can be reconstructed exactly using the filtered back-projection formula (see, for instance, \cite{Buzug08,Feeman15,KakSlaney88}). This formula gives interactions between the Radon transform, the Fourier transform and the back-projection transform. A description of the filtered back-projection formula along \cite[Chapter 3]{KakSlaney88} is provided below.

In the Cartesian system with $x,y$-axes consider a unit vector $(\cos\theta,\sin\theta)$. Then the line perpendicular to this vector with the distance $t$ to the origin can be expressed as ${{\ell }_{t,\theta }}$: $x\cos \theta +y\sin \theta =t$. Assume the object is represented by a two variable function $\mu(x,y)$, which denotes the attenuation coefficient in X-ray CT applications. Then, the $\theta $ -view projection along the line ${{\ell }_{t,\theta }}$ can be expressed as
$$
P(t,\theta )=\int\limits_{-\infty }^{\infty }{\int\limits_{-\infty }^{\infty }{\mu(x,y)\delta (x\cos \theta +y\sin \theta -t)\d x\d y}},
$$
where $\delta$ denotes the Dirac delta-function. The function $P(t,\theta )$ is known as the Radon transform of $\mu(x,y)$. A projection is formed by combining a set of line integrals. The simplest projection is a collection of parallel ray integrals as is given by $P(t,\theta )$ for a constant $\theta $. This is known as a parallel beam projection. It should be noted that there are fan-beam in 2D and cone-beam in 3D projections \cite{Buzug08,Feeman15,KakSlaney88}.
	
The problem of CT is to reconstruct the function $\mu(x,y)$ from its  projections $P(t,\theta )$. There are analytic and iterative methods for CT reconstruction. One of the widely used analytic methods of CT reconstruction is the filtered back-projection method. It can be modeled by
\begin{equation}\label{eq4.8_1}
\mu(x,y)=\int\limits_{0}^{\pi }{\int\limits_{-\infty }^{\infty }{S(\omega ,\theta )\left| \omega  \right|{{e}^{2\pi \i\omega
(x\cos \theta +y\sin \theta )}}\d\omega \d\theta ,}}
\end{equation}
where
\begin{equation}\label{eq4.8_2}
S(\omega ,\theta )=\int\limits_{-\infty }^{\infty }{P(t,\theta ){{e}^{-2\pi \i\omega t}}\d t}		
\end{equation}
is the 1D Fourier transform of $P(t,\theta )$. The inner integral of (\ref{eq4.8_1}) can be regarded as a 1D inverse Fourier transform of the product $S(\omega ,\theta )\left| \omega  \right|$, i.e.,
\begin{equation}\label{eq4.8_3}
Q(t,\theta )=\int\limits_{-\infty }^{\infty }{S(\omega ,\theta )|\omega|{{e}^{2\pi \i\omega t}}\d\omega}
\end{equation}		
which represents a projection filtered by a 1D filter whose frequency representation is $|\omega|$. The outer integral performs back-projection. Therefore, the filtered back-projection consists of two steps: filtration and then back-projection.

Thus, in (\ref{eq4.8_1})-(\ref{eq4.8_3}) the Fourier transforms play the main role. But in practice, due to the fact that we have discrete values of the Radon transform, we have to approximately calculate the Fourier transforms in the filtered back-projection.

Here, in the examples of two functions, we first show that the optimal quadrature formula (\ref{eq3.2})
can be used for approximation of the Fourier transforms of these functions and reconstruction of them.
Then, the optimal quadrature formula (\ref{eq3.2}) is applied for an approximate reconstruction of
the $512\times 512$ size Shepp-Logan phantom from its Radon transform.

Suppose that we are given the values $f(h\beta+a)$, where $\beta=0,1,...,N$, $h=\frac{b-a}{N}$ for $N=2,3,...$.
The values of the function $f$ are assumed to be zero outside the interval $[a,b]$.
Fourier transform  of $f$ (\ref{eq4.4}) is then approximated
by the optimal quadrature formula (\ref{eq3.2}) using $f(h\beta+a)$, $\beta=0,1,...,N$, as follows:
\begin{eqnarray*}
F(\omega)&=&\int\limits_{-\infty}^{\infty}\E^{-2\pi \i \omega x}f(x)\d x\nonumber \\
&=&\int\limits_a^b\E^{-2\pi \i \omega x}f(x)\d x\cong \sum\limits_{\beta=0}^N{C}_{\beta,-\omega}[a,b]f(h\beta+a).
\end{eqnarray*}
Since the coefficients ${C}_{\beta,-\omega}[a,b]$, $\beta=0,1,...,N$, defined by (\ref{eq3.3})-(\ref{eq3.4})
are continuous functions of the variable $\omega$, the following approximation for the Fourier transform
is obtained from the last relations
\begin{equation}\label{eq4.9}
F(\omega)\cong F_{\mathrm{app}}(\omega),
\end{equation}
where
$$
F_{\mathrm{app}}(\omega)=\sum\limits_{\beta=0}^N{C}_{\beta,-\omega}[a,b]f(h\beta+a) \mbox{ for  } \omega\in \mathbb{R}.
$$
Now for the approximate reconstruction of $f$ in the interval $[a,b]$, we approximate the inverse Fourier transform (\ref{eq4.5}) using
the values $F_{\mathrm{app}}(\tau\gamma+a)$, $\gamma=0,1,...,M,$ $\tau=\frac{b-a}{M}$, $M=2,3,...,$ of the function $F_{\mathrm{app}}(\omega)$
in the interval $[a,b]$ for $\omega$ and truncate the integral outside $[a,b]$ as follows:
\begin{eqnarray*}
f(x)=\int\limits_{-\infty}^{\infty}\E^{2\pi \i \omega x}F(\omega)\d\omega
&\cong&\int\limits_{-\infty}^{\infty}\E^{2\pi \i \omega x}F_{\mathrm{app}}(\omega)\d\omega\\
&\cong&\int\limits_a^b\E^{2\pi \i \omega x}F_{\mathrm{app}}(\omega)\d\omega\\
&\cong&\sum\limits_{\gamma=0}^M{C}_{\gamma,x}[a,b]F_{\mathrm{app}}(\tau\gamma+a),
\end{eqnarray*}
where ${C}_{\gamma,x}[a,b]$ are optimal coefficients defined by (\ref{eq3.3})-(\ref{eq3.4}).
Hence, due to the continuity of the coefficients ${C}_{\gamma,x}[a,b]$ for any
$x\in \mathbb{R}$, we obtain an approximation
\begin{equation}\label{eq4.10}
f(x)\cong f_{\mathrm{app}}(x),
\end{equation}
where
$$
f_{\mathrm{app}}(x)=\sum\limits_{\gamma=0}^M{C}_{\gamma,x}[a,b]F_{\mathrm{app}}(\tau\gamma+a)\mbox{ for } x\in \mathbb{R}.
$$
Thus, the function $f(x)$ can be approximately reconstructed, especially in the interval $[a,b]$.
Therefore, formulas (\ref{eq4.9}) and (\ref{eq4.10}) can be used for approximate reconstruction a function from
a set of its values.

\begin{figure}
  \includegraphics[width=\textwidth]{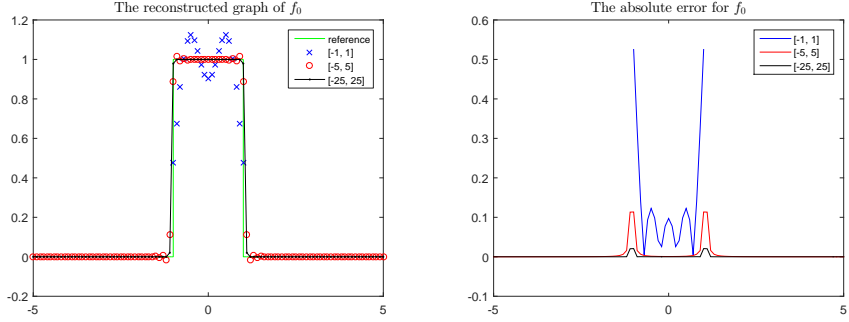}\\
  \caption{The reconstructed graphs and corresponding absolute errors of the function $f_0$ for the intervals $[-1,1]$,
  $[-5,5]$ and $[-25,25]$ using approximation formulas (\ref{eq4.9}) and (\ref{eq4.10}) with steps $h=0.1$ and $\tau=0.01$. }\label{Fig4.4}
\end{figure}

\begin{figure}
\includegraphics[width=\textwidth]{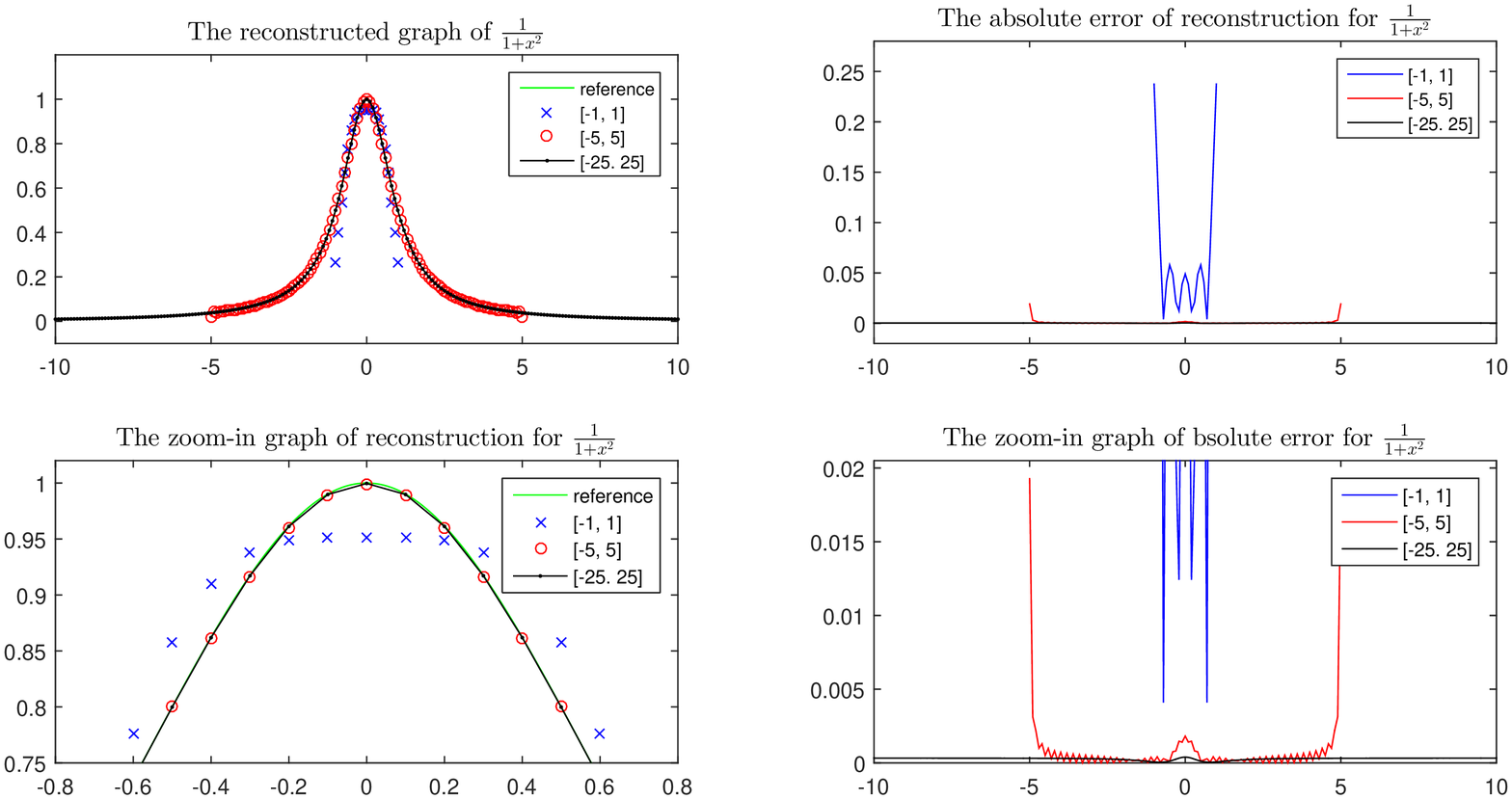}\\
\caption{The reconstructed graphs and corresponding absolute errors of the function $\phi$ for the intervals
$[-1,1]$, $[-5,5]$ and $[-25,25]$ using approximation formulas (\ref{eq4.9}) and (\ref{eq4.10}) with $h=0.1$ and $\tau=0.01$, respectively. }\label{Fig4.5}
\end{figure}

We now demonstrate this in the example of two piecewise continuous functions
$$
f_0(x)=
\left\{
\begin{array}{ll}
1& \mbox{ for } \ x\in [-1,1],\\
0& \mbox{ otherwise}
\end{array}
\right.
$$
and
$$
\phi(x)=
\left\{
\begin{array}{ll}
\frac{1}{1+x^2}& \mbox{ for } \ x\in [a,b],\\
0& \mbox{ otherwise}.
\end{array}
\right.
$$

For reconstruction of these functions, we use the approximation formulas (\ref{eq4.9}) and (\ref{eq4.10}).
In numerical calculations we take the intervals $[-1,1]$, $[-5,5]$ and $[-25,25]$ as an interval of integration $[a,b]$ in
(\ref{eq4.9}) and (\ref{eq4.10}) with $h=0.1$ (for $x$) and $\tau=0.01$ (for $\omega$), respectively.
Then we obtain the reconstructed graphs and graphs of corresponding absolute errors of the functions $f_0$ and $\phi$
for the intervals $[-1,1]$, $[-5,5]$ and $[-25,25]$, shown in Figures \ref{Fig4.4} and \ref{Fig4.5}, respectively.
These graphs show that we get more accurate reconstructions by taking wider intervals.
For the function $f_0$, there are the maximum errors around the jump points $-1$ and $1$ in Figure \ref{Fig4.4}, while for the function $\phi$ in Figure \ref{Fig4.5} the maximum errors are found around the end points of the integration intervals $[-1,1]$, $[-5,5]$ and $[-25,25]$.
Thus, the optimal quadrature formula (\ref{eq3.2}) constructed using coefficients (\ref{eq3.3}) and (\ref{eq3.4})
can be effectively applied to the approximation of Fourier integrals.

Finally, we provide the results of applying the optimal quadrature formula (\ref{eq3.2})
for approximate reconstruction of the $512\times 512$ size Shepp-Logan phantom from its Radon transform.

\begin{figure}
\centering
\includegraphics[width=\textwidth]{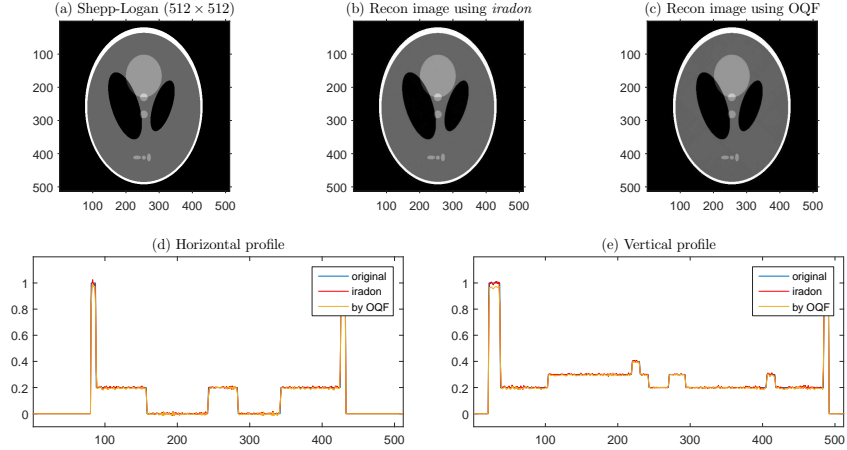}
\caption{Comparison result of CT image reconstruction: (a) Shepp-Logan phantom, (b) reconstructed CT image using MATLAB built-in function \textit{iradon}, (c) reconstructed CT image using optimal quadrature formula (OQF), (d) profile of horizontal centerline, (e) profile of vertical centerline.\label{recon_results}}\end{figure}

We generate the sinogram using half rotation sampling with sampling angle~$0.5^{\circ}$.
We compare the result of CT image reconstruction using the optimal quadrature formula for Fourier integrals with the result of  \textit{iradon}, a built-in function of MATLAB R2019a, which uses \textit{fft} and \textit{ifft} for Fourier integrals.
For the image quality analysis, we compare maximum error~($E_{\max}$), mean squared error (MSE), and the peak signal-to-noise ratio~(PSNR):
\begin{eqnarray*}
E_{\max}(I)&=&\max_{i,j}|I(i,j)-I_{ref}(i,j)|, \\
\mbox{MSE}(I)&=&\frac{1}{mn}\sum_{i=1}^{m}\sum_{j=1}^{n}|I(i,j)-I_{ref}(i,j)|^2\,, \\
\mbox{PSNR}(I)&=&10\log_{10}   \left( \frac{I_{\max}^2}{\mbox{MSE}(I)} \right)\,,
\end{eqnarray*}
where
$I_{\max}$ is the maximum pixel value of the image $I$.
For $I_{ref}$, we adopt a Shepp-Logan phantom (Figure~\ref{recon_results}(a)).

\begin{center}
\begin{table}[]
\caption{Quantitative image analysis}
\begin{tabular}{|l||c|c|c|}
\hline
& \multicolumn{1}{c|}{\begin{tabular}[c]{@{}c@{}}$E_{\max}$\\ (inner part)\end{tabular}} & \multicolumn{1}{c|}{\begin{tabular}[c]{@{}c@{}}MSE\\ (inner part)\end{tabular}} & \multicolumn{1}{c|}{\begin{tabular}[c]{@{}c@{}}PSNR \\ (inner part)\end{tabular}} \\ \hline\hline
\begin{tabular}[c]{@{}l@{}}Reconstruction result using\\ MATLAB built-in function \\ \textit{iradon} \end{tabular}     & \begin{tabular}[c]{@{}c@{}} 0.3472 \\ (0.2812)    \end{tabular}                 & \begin{tabular}[c]{@{}c@{}}9.2103e-04 \\ (1.9613e-04)\end{tabular}              & \begin{tabular}[c]{@{}c@{}}30.357 \\ (37.075)\end{tabular}                      \\ \hline
\begin{tabular}[c]{@{}l@{}}Reconstruction result using\\ the proposed optimal \\ quadrature formula\end{tabular} & \begin{tabular}[c]{@{}c@{}}0.3895 \\ (0.2689)\end{tabular}                      & \begin{tabular}[c]{@{}c@{}}10.8548e-04 \\ (1.9265e-04)\end{tabular}             & \begin{tabular}[c]{@{}c@{}}29.644 \\ (37.152)\end{tabular}                      \\ \hline
\end{tabular}
\label{imageAnalysis}
\end{table}
\end{center}

As shown in Figure~\ref{recon_results}(c), the CT image reconstruction algorithm using optimal quadrature formula produces a clear reconstruction image which has the same structures with the original phantom.
It also has almost the same appearance with the result of \textit{iradon} as shown in Figure~\ref{recon_results}(b).
From Figure~\ref{recon_results}(d) and~(e), we see that the results using the optimal quadrature formula and \textit{iradon} are almost the same except the outer ring.
Table~\ref{imageAnalysis} shows $E_{\max}$, MSE, and PSNR for two reconstruction results.
The numbers written without parentheses are measured errors in the whole image domain and those in parentheses are measured ones only inside the outer ring.
The reconstruction result by using \textit{iradon} seems better than the result by using our optimal quadrature formula if we consider the whole image domain for the error measurement.
However, in most CT applications, we are interested in interior structures of an object rather than its outer part, hence we exclude the outer ring, then the opposite holds true.
Note that unlike \textit{iradon} which uses various optimized image processing techniques as a MATLAB built-in function, no image processing technique has been applied to our image reconstruction process.
We expect that the performance can be improved further if we use the higher order optimal quadrature formula, which is our next research topic.

\section*{Acknowledgements}
The work has been done while A.R. Hayotov was visiting Department of Mathematical
Sciences at KAIST, Daejeon, Republic of Korea. A.R. Hayotov's work was supported by the 'Korea Foundation for Advanced Studies'/'Chey Institute for Advanced Studies' International Scholar Exchange Fellowship for academic year of 2018-2019. C.-O. Lee's work was supported by NRF grant funded by MSIT (NRF-2017R1A2B4011627).


\newpage

\begin{thebibliography}{99}
\bibitem{AvdMal89} V.A. Avdeenko and A.A. Malyukov, A quadrature formula for the Fourier integral based on the use of a cubic spline,
(Russian) USSR Computational Mathematics and Mathematical Physics,  29 (1989) 783-786.

\bibitem{IBab}  I. Babu\v{s}ka, Optimal quadrature formulas (Russian), Dokladi
Akad. Nauk SSSR,  149 (1963), 227--229.

\bibitem{BabVitPrag69} I. Babu\v{s}ka, E. Vitasek, and M. Prager, Numerical processes in differential equations. Wiley, New York, 1966.

\bibitem{BakhVas68} N.S. Bakhvalov and L.G. Vasil'eva, Evaluation of the integrals of oscillating functions by interpolation at nodes
of Gaussian quadratures (Russian) USSR Computational Mathematics and Mathematical Physics, 8 (1968), 241-249.

\bibitem{BolHayShad16} N.D. Boltaev, A.R. Hayotov, and Kh.M. Shadimetov, Construction of optimal quadrature formula for numerical calculation
of Fourier coefficients in Sobolev space $L_2^{(1)}$, American Journal of Numerical Analysis, 4 (2016), 1-7.

\bibitem{BolHayShad17} N.D. Boltaev, A.R. Hayotov and Kh.M. Shadimetov, Construction of optimal quadrature formulas for Fourier coefficients in Sobolev space $L_2^{(m)}(0,1)$, Numerical Algorithms, 74 (2017),  307-336.

\bibitem{BolHayMilShad17} N.D. Boltaev, A.R. Hayotov, G.V. Milovanovi\'{c}, and Kh.M. Shadimetov, Optimal quadrature formulas for Fourier coefficients
in $W_2^{(m,m-1)}$ space, Journal of applied analysis and  computation, 7 (2017), 1233-1266.

\bibitem{Buzug08} T.M. Buzug, Computed tomography, From photon statistics to modern cone-beam CT, Springer, Berlin, 2008.

\bibitem{CatCom}  T. Catina\c{s} and Gh. Coman, Optimal quadrature formulas based on the $\phi$-function method, Stud. Univ. Babe\c{s}-Bolyai Math., 51
(2006), 49-64.

\bibitem{Feeman15} T.G. Feeman, The mathematics of medical imaging, A Beginner's guide, Second edition, Springer, Switzerland, 2015.

\bibitem{Filon28} L.N.G. Filon, On a quadrature formula for trigonometric integrals, Proc. Roy. Soc. Edinburgh, 49 (1928), 38-47.

\bibitem{GhOs}  A. Ghizzetti and A. Ossicini,  Quadrature Formulae,
 Akademie Verlag, Berlin, 1970.


\bibitem{IserNor05} A. Iserles and S.P. N{\o}rsett. Efficient quadrature of highly oscillatory integrals using derivatives, Proc. R. Soc. A,
461 (2005), 1383--1399.

\bibitem{KakSlaney88} A.C. Kak and M. Slaney, Principles of Computerized Tomographic imaging, IEEE Press, New York, 1988.

\bibitem{FLan}  F. Lanzara, On optimal quadrature formulae, J. Ineq. Appl., 5 (2000),  201-225.

\bibitem{Mil98} G.V. Milovanovi\'c, Numerical calculation of integrals involving oscillatory and singular kernels and
some applications of quadratures, Computers Math. Applic., 36 (1998), 19-39.

\bibitem{MilStan14} G.V. Milovanovi\'c and M.P. Stani\'c, Numerical integration of highly oscillating functions, In:
Analytic Number Theory, Approximation Theory, and Special Functions/ G.V. Milovanovi\'c amd M.Th. Rassias (Eds.), Springer, New York,
2014, 613-649.

\bibitem{NovUllWoz15} E. Novak, M. Ullrich, and H. Wo\'zniakowski, Complexity of oscillatory integration for univariate
Sobolev space, Journal of Complexity, 31 (2015), 15-41.


\bibitem{Olver08} S. Olver, Numerical approximation of highly oscillatory integrals, PhD dissertation, University of Cambridge, 2008.

\bibitem{Sard}  A. Sard, Best approximate integration formulas; best approximation formulas, Amer. J. Math., 71 (1949), 80-91.

\bibitem{Shad85}    Kh.M. Shadimetov, The discrete analogue of the
differential operator $\d^{2m}/\d x^{2m}$ and its construction,
Questions of Computations and Applied Mathematics. Tashkent, no. 79
(1985), 22-35. ArXiv:1001.0556.v1 [math.NA] Jan. 2010.

\bibitem{Shad99}
 Kh.M. Shadimetov, Weight optimal cubature formulas in Sobolev's periodic space. (Russian)
Siberian J. Numer. Math. -Novosibirsk, 2 (1999), 185-196


\bibitem{ShadHay11}  Kh.M. Shadimetov and A.R. Hayotov, Optimal quadrature formulas with positive coefficients in $L_2^{(m)}(0,1)$ space,
J. Comput. Appl. Math., 235 (2011), 1114-1128.

\bibitem{ShadHayAkhm15} Kh.M. Shadimetov, A.R. Hayotov, and D.M. Akhmedov, Optimal quadrature formulas for Cauchy type singular integrals
in Sobolev space, Applied Mathematics and Computation, 263 (2015), 302-314.

\bibitem{Sobolev06}  S.L. Sobolev, The coefficients of optimal quadrature formulas, Selected Works of S.L. Sobolev,  Springer US,  2006, 561-566.

\bibitem{Sobolev74} S.L. Sobolev, Introduction to the  theory of cubature formulas (Russian), Nauka, Moscow,  1974.

\bibitem{SobVas} S.L. Sobolev and  V.L. Vaskevich, The theory of cubature formulas, Kluwer Academic Publishers Group, Dordrecht, 1997.

\bibitem{XuMilXiang} Z. Xu, G.V. Milovanovi\'c, and S. Xiang, Efficient computation of highly oscillatory integrals with
Henkel kernel, Appl. Math. and Comp., 261 (2015), 312-322.

\bibitem{ZhangNovak19} S. Zhang and E. Novak, Optimal quadrature formulas for the Sobolev space $H^1$, Journal of Scientific Computing,
78 (2019), 274-289.
\end{thebibliography}
\end{document}